\documentclass[11pt]{article}

\makeindex
\usepackage[all]{xy}
\usepackage{graphicx}

\usepackage{amsmath}
\usepackage{amsfonts}
\usepackage{amssymb}
\usepackage{amsthm}
\newcommand{\isom}{\cong}

\newcommand{\Z}{{\bf{Z}}}

\newcommand{\Q}{{\bf{Q}}}
\newcommand{\Qbar}{{\overline{\Q}}}

\newcommand{\C}{{\bf{C}}}
\newcommand{\F}{{\bf{F}}}
\newcommand{\T}{{\bf{T}}}

\newcommand{\m}{{\mathfrak{m}}}
\newcommand{\OO}{\mathcal{O}}
\newcommand{\NN}{\mathcal{N}}

\newcommand{\noi}{{\noindent}}

\newcommand{\Mid}{|} %{\mid\!\!}

\newcommand{\miD}{|}%{\!\!\mid}

\newcommand{\divs}{\!\mid\!}
\newcommand{\ndiv}{\!\nmid\!}
\newcommand{\tensor}{\otimes}

\newcommand{\ra}{{\rightarrow}}

 1
\DeclareFontEncoding{OT2}{}{} % to enable usage of cyrillic fonts
  \newcommand{\textcyr}[1]{%
    {\fontencoding{OT2}\fontfamily{wncyr}\fontseries{m}\fontshape{n}%
     \selectfont #1}}
\newcommand{\Sha}{{\mbox{\textcyr{Sh}}}} %{\hbox{\cyr X}}

\newcommand{\Adual}{\widehat{A}}
\newcommand{\Agdual}{A_g^{\vee}}
\newcommand{\Ahdual}{A_h^{\vee}}
\newcommand{\Afdual}{A_f^{\vee}}
\newcommand{\Edual}{E^{\vee}}
\newcommand{\Fdual}{F^{\vee}}

\newcommand{\pp}{{\mathfrak{p}}}
\newcommand{\qq}{{\mathfrak{q}}}

\newcommand{\annT}{{{\rm Ann_\T}}}

\newcommand{\tphat}{{\widehat{\T P}}}

\newcommand{\comment}[1]{}
\newcommand{\marginalfootnote}[1]{%
   \footnote{#1}\marginpar{\hfill {\sf\thefootnote}}%
}

\newcommand{\edit}[1]{\marginalfootnote{#1}}

\newtheorem{lem}{Lemma}[section]
\newtheorem{cor}[lem]{Corollary}
\newtheorem{prop}[lem]{Proposition}
\newtheorem{conj}[lem]{Conjecture}
\newtheorem{thm}[lem]{Theorem}

\theoremstyle{definition}

\newcommand{\thetitle}
{Visibility and the Birch and Swinnerton-Dyer conjecture
for analytic rank one}

\begin{document}
%\ssp
\parindent=2em

\title{\thetitle}
\author{Amod Agashe
\footnote{This material is based upon work supported by the National Science 
Foundation under Grant No. 0603668.}}
\maketitle

%\include{abstract}
%abstract.tex
\begin{abstract}
Let $E$ be an optimal elliptic curve over~$\Q$ of conductor~$N$
having analytic rank one, i.e., 
such that the $L$-function~$L_E(s)$ of~$E$ vanishes to order one at~$s=1$.
Let $K$ be a quadratic imaginary field in which all 
the primes dividing~$N$ split and such that 
the $L$-function of~$E$ over~$K$ vanishes to order one at~$s=1$. 
Suppose there is another optimal elliptic curve over~$\Q$ of the same conductor~$N$
whose  Mordell-Weil rank is greater than one and
whose associated newform is congruent to the newform associated to~$E$
modulo an integer~$r$. The theory of visibility 
then shows that under certain additional hypotheses,
$r$ divides the order of the Shafarevich-Tate group of~$E$ 
over~$K$. We show that under somewhat similar hypotheses, 
$r$ also divides the Birch and Swinnerton-Dyer
{\em conjectural} order of the Shafarevich-Tate group of~$E$ over~$K$, 
which provides new theoretical evidence for
the second part of 
the Birch and Swinnerton-Dyer conjecture in the analytic rank one case.
\end{abstract}

%Warning: This is a very rough draft -- I have mainly just
%jotted down some ideas. At several places,
%one may need minor hypotheses that have been skipped for simplicity.

\section{Introduction and results}

Mazur introduced the notion of visibility in order to better 
understand geometrically the elements of the 
Shafarevich-Tate group of an abelian variety.
The correspoding theory, which we call the theory of visibility, 
can often be used to show the existence of 
non-trivial elements of the Shafarevich-Tate group
of abelian varieties and motives (e.g., see~\cite{cremona-mazur},
\cite{agst:vis}, \cite{dsw}). 
The second part of the
Birch and Swinnerton-Dyer conjecture gives a 
formula for the order of the Shafarevich-Tate group,
and one might wonder how much of this conjectural order,
when non-trivial, 
can be explained by the theory of visibility.
This issue has been investigated computationally
(e.g., see~\cite{cremona-mazur}, \cite{agst:bsd})
and theoretically (e.g.,\cite{agmer}) when the concerned
abelian variety has analytic rank zero,
but not for any higher analytic rank 
(where, by the {\em analytic rank} of
an abelian variety, we mean
the order of vanishing of the $L$-function of the abelian variety
at $s=1$).
%, assuming the analytic continuation of the $L$-function, which
%will always hold in our case).
%(as far as we know). 
In this article, we take the first step
in investigating the issue for higher analytic rank, by
showing that for elliptic curves of analytic rank one,
in certain situations
where the theory of visibility implies that
the actual Shafarevich-Tate group is non-trivial,
the Birch and Swinnerton-Dyer conjectural order of
the Shafarevich-Tate group is non-trivial  as well.
%relating the Birch and Swinnerton-Dyer conjectural order of
%the Shafarevich-Tate group to its actual order using
%the theory of visibility for elliptic curves of analytic rank one.
This provides new theoretical evidence for
the second part of
the Birch and Swinnerton-Dyer conjecture in the analytic rank one case,
and links the theory of visibility to the 
Birch and Swinnerton-Dyer conjecture for the first time
for any elliptic curve of analytic rank greater than zero.
%and gives the hope that the theory of visibility may give useful information
%even when the analytic rank is greater than zero. 

We now state our results more precisely.
Let $N$ be a positive integer. Let
%, which we might have to assume
%is squarefree in several situations below, 
$X_0(N)$ be the modular curve over~$\Q$
associated to~$\Gamma_0(N)$, and
let $J=J_0(N)$ denote the Jacobian of~$X_0(N)$, which is
an abelian variety over~$\Q$. 
%If $f$ is an eigenform of weight~$2$ on~$\Gamma_0(N)$,
%we call the order of vanishing of the $L$-function $L(f,s)$ at~$s=1$
%the {\em analytic rank of~$g$}.
%If $f$ be a newform of weight~$2$ on~$\Gamma_0(N)$ 
%whose analytic rank is one.
Let $\T$ denote the Hecke algebra, which is 
the subring of endomorphisms of~$J_0(N)$
generated by the Hecke operators (usually denoted~$T_\ell$
for $\ell \ndiv N$ and $U_p$ for $p\divs N$). 
If $f$ is a newform of weight~$2$ on~$\Gamma_0(N)$, then 
let $I_f = \annT f$ and let $A_f$ denote the associated {\em newform
quotient} $J/I_f J$, which is
an abelian variety over~$\Q$. 
Let $\pi$ denote the quotient map $J \ra J/I_f J = A_f$.
%Suppose that $A_f$ has analytic rank one.
By the {\em analytic rank} of~$f$, we mean the
order of vanishing at $s=1$ of~$L(f,s)$. The analytic rank of~$A_f$
is then the analytic rank of~$f$ times the dimension of~$A_f$.
Now suppose that the newform~$f$ has integer Fourier coefficients.
Then $A_f$ is an elliptic curve, and we denote it by~$E$ instead.
Since $E$ has dimension one, it has analytic rank one.
%, i.e., that the $L$-function~$L_E(s)$ of~$E$ vanishes to order one at~$s=1$
%(equivalently, that the $L$-function $L(f,s)$ of~$f$
%vanishes to order one at~$s=1$).

%Suppose $D \neq -3, -4$ is a negative fundamental discriminant,
%and let $K = \Q(\sqrt{D})$ be the associated quadratic imaginary field.
Let $K$ be a quadratic imaginary field of discriminant not equal
to~$-3$ or~$-4$, and  such that
all primes dividing~$N$
split in~$K$. Choose an ideal~$\NN$ of the ring of integers~$\OO_K$
of~$K$ such that $\OO_K/\NN \isom \Z/N \Z$. Then the complex tori
$\C/\OO_K$ and~$\C/\NN^{-1}$ define elliptic curves related by 
a cyclic $N$-isogeny, and thus give a complex valued point~$x$ of~$X_0(N)$.
This point, called a Heegner point, is defined over the
Hilbert class field~$H$ of~$K$.
%Let $x$ be a Heegner point of discriminant~$D$ on~$X$ (as in~\cite[\S~I.3]{gross-zagier}). 
Let $P \in J(K)$ be the 
class of the divisor 
$\sum_{\sigma \in {\rm Gal}(H/K)} ((x) - (\infty))^\sigma$, where
$H$ is the Hilbert class field of~$K$. 
\comment{
Suppose $D \neq -3, -4$ is a negative fundamental discriminant,
and let $K = \Q(\sqrt{D})$. Let $x$ be a Heegner point of
discriminant~$D$ on~$X$ (as in~\cite[\S~I.3]{gross-zagier}). 
Let $P \in J(K)$ be the 
class of the divisor 
$\sum_{\sigma \in {\rm Gal}(H/K)} ((x) - (\infty))^\sigma$, where
$H$ is the Hilbert class field of~$K$. 
}

By~\cite{waldspurger:comp2}, we may choose K so 
that $L(E/K,s)$ vanishes to order one at $s=1$.
%twisted L-value is nonzero. So L'(E/K,1) is non-zero,
Hence, by~\cite[{\S}V.2:(2.1)]{gross-zagier}, $\pi(P)$ has 
infinite order, and by work of Kolyvagin, $E(K)$ has rank one and
the order of the Shafarevich-Tate group~$\Sha(E/K)$ of~$E$ over~$K$ 
is finite
(e.g., see~\cite[Thm.~A]{kolyvagin:euler} or~\cite[Thm.~1.3]{gross:kolyvagin}).
In particular, the index $[E(K):\Z\pi(P)]$
is finite.
%is finite
By~\cite[{\S}V.2:(2.2)]{gross-zagier}
(or see~\cite[Conj.~1.2]{gross:kolyvagin}), the second part of
the Birch and Swinnerton-Dyer (BSD) conjecture becomes:
\begin{conj}[Birch and Swinnerton-Dyer, Gross-Zagier] \label{conj:bsd}
\begin{eqnarray} \label{gzformula}
|E(K)/ \Z \pi(P)| 
\stackrel{?}{=} c_{\scriptscriptstyle E} \cdot \prod_{p | N} c_p(E) \cdot 
\Mid \Sha(E/K) \miD^{1/2},
\end{eqnarray}
where $c_{\scriptscriptstyle E}$ is the Manin constant of~$E$,
$c_p(E)$ denotes the arithmetic component group of~$E$ at the prime~$p$,
and the question mark above the equality sign emphasizes that this
equality is conjectural.
\end{conj}

Note that the Manin constant~$c_{\scriptscriptstyle E}$
is conjectured to be one, and
one knows that if $p$ is a prime such that $p^2 \nmid 4 N$, 
then $p$ does not divide~$c_{\scriptscriptstyle E}$
(by~\cite[Cor.~4.1]{mazur:rational} 
and~\cite[Thm.~A]{abbes-ullmo}).

Now suppose that $f$ is congruent to another newform~$g$ with integer
Fourier coefficients, 
whose associated elliptic curve~$F$ has Mordell-Weil
rank over~$\Q$ bigger than one. Then the theory of visibility 
(e.g., as in~\cite{cremona-mazur})
often shows
that such congruences divide the order of~$\Sha(E/K)$,
as we now indicate.

We say that a maximal ideal~$\m$ of~$\T$
satisfies {\em multiplicity one} if
$J_0(N)[\m]$ is two dimensional over~$\T/\m$.
This is known to hold in several situations, in particular
when the following conditions hold simultaneously
(e.g., see~\cite[Thm.~2.1(ii)]{wiles:fermat} along with~\cite{ars:moddeg}):
$p \neq 2$, $p^2 \ndiv N$, 
the canonical semi-simple representation~$\rho_\m$
associated to~$\m$ (see, e.g., \cite[Prop.~5.1]{ribet:modreps}
for the defintion of~$\rho_\m$)
is irreducible, and $\m$ arises as a pullback
from $\T/I_h$ for some newform~$h$.

\begin{lem} \label{lem:containment}
Let $r$ be an integer such that:\\
(a) $f$ and~$g$  are congruent modulo~$r$ and 
for every prime~$p$ that divides~$r$, $f$ and~$g$ are not
congruent modulo a power of~$p$ greater than~$p^{{\rm ord}_p r}$. \\
(b) if  $\m$ is a maximal ideal of~$\T$ such that the residue characteristic
of~$\m$ divides~$r$ and $\m$ is in the support of~$J_0(N)[I_f]$
or~$J_0(N)[I_g]$, then $\m$ satisfies multiplicity one.\\
Then $\Edual[r] = \Fdual[r]$, and both are direct summands
of~$\Edual \cap \Fdual$ as~${\rm Gal}(\Qbar/\Q)$-modules.
\end{lem}
\begin{proof}
By~\cite[Cor.~2.5]{emerton:optimal},
if $\m$ satisfies multiplicity one 
and $I$ is any saturated ideal of~$\T$, then 
the $\m$-adic completion of the group of connected components of~$J_0(N)[I]$
is trivial. 
If $L \ra M$ is a homomorphism
of two $\T$-modules,  
then we say that $L = M$ {\em away from} a given 
set of maximal ideals
if the induced map on the $\m$-adic completions is an isomorphism
for all maximal ideals~$\m$ that are not in the prescribed set.
Thus, the inclusions $\Edual \subseteq J_0(N)[I_f]$
and $\Fdual \subseteq J_0(N)[I_g]$ are equalities
away from maximal ideals that do not satisfy multiplicity one.
Hence $\Edual \cap \Fdual \subseteq J_0(N)[I_f + I_g]$ is an
equality away from maximal ideals that do not satisfy multiplicity one.
But $\Edual \cap \Fdual \subseteq \Edual[I_f + I_g] 
\subseteq J_0(N)[I_f + I_g]$. 
Hence if $m$ denotes the largest integer
such that $f$ and~$g$ are congruent modulo~$m$, 
then $\Edual \cap \Fdual \subseteq \Edual[I_f + I_g] = \Edual[m]$
is an equality away from maximal ideals that do not satisfy multiplicity one.
Similarly 
$\Edual \cap \Fdual \subseteq \Fdual[I_f + I_g] = \Fdual[m]$
is an equality away from maximal ideals that do not satisfy multiplicity one.
%Let $r'$ be the maximal divisor of~$r$ such that 
%for all maximal ideals~$\m$ of~$\T$ whose residue characteristic
%divides~$r'$ satisfy multiplicity one.
From conditions~(a) and~(b) on~$r$, it 
follows that for every prime~$p$ that divides~$r$,
$(E \cap F)[p^\infty] = \Edual[p^\infty] = \Fdual[p^\infty]$.
Thus
$\Edual[r]$ and~$\Fdual[r]$ are identical and
are direct summands 
of~$\Edual \cap \Fdual$ as~${\rm Gal}(\Qbar/\Q)$-modules.
\end{proof}

\begin{prop} \label{prop:agst}
(i) Let $r$ be the largest integer such that 
%the newforms associated
%to~$E$ and~$F$ 
$f$ and~$g$ are congruent modulo~$r$ and
all maximal ideals~$\m$ of~$\T$ such that the residue characteristic
of~$\m$ divides~$r$ and $\m$ is in the support of~$J_0(N)[I_f]$
or~$J_0(N)[I_g]$
satisfy multiplicity one.
Suppose that $r$ is coprime to 
$$N \cdot |(J_0(N)/\Fdual)(K)_{\rm tor}|
\cdot |F(K)_{\rm tor}| 
\cdot \prod_{\scriptscriptstyle{{p\mid N}}} c_p(F).$$
Then $r$ divides~$\prod_{\scriptscriptstyle{{p\mid N}}} c_p(E) \cdot |\Sha(E/K)|^{1/2}$,
which in turn divides the right hand side of
the Birch and Swinnerton-Dyer conjectural formula~(\ref{gzformula}).\\
(ii) Suppose that $f$ is congruent to~$g$ modulo an odd prime~$q$
such that $E[q]$ and~$F[q]$ are irreducible and 
$q$ does not divide 
$$N \cdot |(J_0(N)/\Fdual)(K)_{\rm tor}| \cdot |F(K)_{\rm tor}|.$$
%\cdot \prod_{\scriptscriptstyle{{p\mid N}}} c_p(F).$$
Also, assume that $f$ is not congruent modulo~$q$ to a newform
of a level dividing~$N/p$ for some prime~$p$ that divides~$N$
(for Fourier coefficients of index coprime to~$Nq$),
and either $q \nmid N$ or for all primes $p$ that divide~$N$,
$q \nmid (p-1)$. Then $q$ 
%~$c_{\scriptscriptstyle E} \cdot \prod_{\scriptscriptstyle{{p\mid N}}} c_p(E)
%\cdot 
divides~$\Mid \Sha(E/K) \miD$.
%, the right hand side of
%the Birch and Swinnerton-Dyer conjectural formula~(\ref{gzformula}).
\end{prop}
\begin{proof}
Both results follow essentially from Theorem~3.1 of~\cite{agst:vis}.
For the first part, take $A = \Edual$, $B = \Fdual$, and $n=r$ 
in~\cite[Thm.~3.1]{agst:vis}, and
note that
$\Fdual[r] \subseteq \Edual$ by Lemma~\ref{lem:containment}
and that the rank of~$\Edual(K)$ 
is less than the rank of~$\Fdual(K)$. 
For the second part, take $A = \Edual$, $B = \Fdual$, and $n=q$ 
in~\cite[Thm.~3.1]{agst:vis}, and note that 
the congruence of $f$ and~$g$ modulo~$q$
forces $\Fdual[q] = \Edual[q]$ 
by~\cite[Thm.~5.2]{ribet:modreps} (cf.~\cite[p.~20]{cremona-mazur}),
and that the hypotheses imply that 
$q$ does not divide~$c_p(E)$ or $c_p(F)$ for any prime~$p$
that divides~$N$, as we now indicate.
By~\cite[Prop.~4.2]{emerton:optimal}, if $q$ divides~$c_p(E)$ 
for some prime~$p$ that divides~$N$, then for some maximal ideal~$\m$
of~$\T$ having characteristic~$q$ and containing~$I_f$, 
either $\rho_\m$ is finite or reducible.
Since $E[q]$ is irreducible, this can happen only if 
$\rho_\m$ is finite. But 
this is not possible by~\cite[Thm.~1.1]{ribet:modreps},
in view of the hypothesis that  $f$
is not congruent modulo~$q$ to a newform
of a level dividing~$N/p$ for some prime~$p$ that divides~$N$
(for Fourier coefficients of index coprime to~$Nq$),
and either $q \nmid N$ or for all primes $p$ that divide~$N$,
$q \nmid (p-1)$. Thus $q$ does not divide~$c_p(E)$ for any prime~$p$
that divides~$N$. Similarly,
$q$ does not divide~$c_p(F)$ for any prime~$p$
that divides~$N$,
considering that the hypothesis that $f$
is not congruent modulo~$q$ to a newform
of a level dividing~$N/p$ for some prime~$p$ that divides~$N$
(for Fourier coefficients of index coprime to~$Nq$) applies
to $g$ as well, since
$g$ is congruent to~$f$ modulo~$q$. This finishes the proof
of the proposition.
\end{proof}

One might wonder how often it happens in numerical
data that visibility explains the Shafarevich-Tate group
of an elliptic curve of analytic rank one. 
Since it is difficult to compute the actual 
order of the Shafarevich-Tate group,
we looked at the Birch and Swinnerton-Dyer conjectural
orders in Cremona's online ``Elliptic curve data''~\cite{cremona:online}.
For levels up to~$30000$, we found only
one optimal elliptic curve of Mordell-Weil rank one for which
the conjectural order of the Shafarevich-Tate group was divisible
by an odd prime: the curve
with label 28042A, for which the
conjectural order of the Shafarevich-Tate group is~$9$. 
At the same level, the curve 28042B has Mordell-Weil rank~$3$ and 
the newforms corresponding to 28042A and 28042B have Fourier
coefficients that are congruent modulo~$3$ for every prime index up to~$100$.
While this is not enough to conclude that the newforms are congruent
modulo~$3$ for all Fourier coefficients 
(cf.~\cite{agashe-stein:schoof-appendix}),
it is  quite likely that this is true  
%are congruent modulo~$3$, 
and that this congruence explains the non-trivial Shafarevich-Tate group,
although we have not checked the details (in particular whether
the hypotheses of Proposition~\ref{prop:agst} are satisfied),
since our goal in this paper is to prove {\em theoretical} results.
It would be interesting to do systematic computations to see
how much of the Birch and Swinnerton-Dyer conjectural order of
the Shafarevich-Tate group is explained by visibility for 
elliptic curves of analytic rank one
(similar to the computations
in~\cite{agst:bsd} for the analytic rank zero case).

In any case, when the theory of visibility does imply
that the Shafarevich-Tate group is non-trivial, this
should be reflected in its Birch and Swinnerton-Dyer 
conjectural order. In our situation, 
since Proposition~\ref{prop:agst} shows that under certain hypotheses,
certain integers of
congruence divide the right side 
of the Birch and Swinnerton-Dyer conjectural formula~(\ref{gzformula}), 
these integers should also divide the left
side of~(\ref{gzformula}),
if the second part of the 
Birch and Swinnerton-Dyer conjecture is true.
We have the following result
in this direction.

\begin{thm} \label{thm:main}
Recall that $E$ and~$F$ are optimal elliptic curves over~$\Q$
associated to newforms~$f$ and~$g$ respectively,
with $E$ of analytic rank one and $F$ having Mordell-Weil rank more
than one. 
Let $r$ be an integer such that:\\
(a) $f$ and~$g$  are congruent modulo~$r$ and 
for every prime~$p$ that divides~$r$, $f$ and~$g$ are not
congruent modulo a power of~$p$ greater than~$p^{{\rm ord}_p r}$. \\
(b) all maximal ideals~$\m$ of~$\T$ whose residue characteristic
divides~$r$ and that are in the support of~$J_0(N)[I_f]$
or of~$J_0(N)[I_g]$ satisfy multiplicity one.\\
(c) $r$ is coprime to the order of the torsion subgroup of
the projection of~$\T P$ in~$J_0(N)/(I_f \cap I_g) J_0(N)$.\\
Then $r$ divides $|E(K)/ \Z \pi(P)|$, which is 
the left hand side of
the Birch and Swinnerton-Dyer conjectural formula~(\ref{gzformula}).
\end{thm}

We will prove this theorem in Section~\ref{sec:proofs}.

\begin{cor} \label{cor:main}
Suppose $f$ and~$g$ are congruent modulo an odd prime~$q$
such that $q^2 \nmid N$, $E[q]$ and~$F[q]$ are irreducible, and
$q$ does not divide $|J_0(N)(K)_{\rm tor}|$. 
Then $q$ divides $|E(K)/ \Z \pi(P)|$.
If moreover, 
$f$ is not congruent modulo~$q$ to a newform
of a level dividing~$N/p$ for some prime~$p$ that divides~$N$
(for Fourier coefficients of index coprime to~$Nq$),
and either $q \nmid N$ or for all primes $p$ that divide~$N$,
$q \nmid (p-1)$, then 
$q$ divides the Birch and Swinnerton-Dyer conjectural
order of~$\Sha(E/K)$.
%If moreover, there exists $K$ such that 
%${rm ord}_q L(E_D,s)/\Omega_{E_D} = 0$, where $D$ is the
%discriminant of~$K$, $E_D$ denotes the twist of~$E$ by
%the Kronecker character associated to~$K$, and $\Omega_{E_D}$
%is the N\'eron period of~$E$, then 
%$q$ divides the Birch and Swinnerton-Dyer conjectural
%order of~$\Sha(E/\Q)$.
\end{cor}
\begin{proof}
Let $r$ denote the highest power of~$q$ modulo which $f$ and~$g$
are congruent. Thus condition (a) on~$r$
in Theorem~\ref{thm:main} is satisfied.
By the discussion just before Lemma~\ref{lem:containment},
the hypotheses that $q$ is odd, $q^2 \nmid N$, and 
$E[q]$ and~$F[q]$ are  irreducible 
imply that $r$ satisfies condition~(b).
The hypothesis that $q$ does not divide 
 $|J_0(N)(K)_{\rm tor}|$ implies that condition~(c) is satisfied.
Hence by Theorem~\ref{thm:main}, $q$ divides $|E(K)/ \Z \pi(P)|$.
%By~\cite[Prop.~4.2]{emerton:optimal}, if $q$ divides~$c_p(E)$ 
%for some prime~$p$ that divides~$N$, then for some maximal ideal~$\m$
%of~$\T$ having characteristic~$q$ and containing~$I_f$, 
%either $\rho_\m$ is finite or reducible.
%Since $E[q]$ is irreducible, this can happen only if 
%$\rho_\m$ is finite. But 
%this is not possible by~\cite[Thm.~1.1]{ribet:modreps},
%in view of the hypothesis that  $f$
%is not congruent modulo~$q$ to a newform
%of a level dividing~$N/p$ for some prime~$p$ that divides~$N$
%(for Fourier coefficients of index coprime to~$Nq$),
%and either $q \nmid N$ or for all primes $p$ that divide~$N$,
%$q \nmid (p-1)$. Thus 
As explained in the proof of Proposition~\ref{prop:agst},
the hypotheses imply that
$q$ does not divide~$c_p(E)$ for any prime~$p$.
%that divides~$N$. 
Also, by~\cite[Cor.~4.1]{mazur:rational},
$q$ does not divide the Manin constant~$c_{\scriptscriptstyle E}$. 
Hence, by~(\ref{gzformula}), 
$q$ divides the Birch and Swinnerton-Dyer conjectural
order of~$\Sha(E/K)$. 
%Finally,
%if $L(E/K,s)$ denotes the $L$-function of~$E$ over~$K$ and
%$\Omega(E/K)$ denotes the N\'eron period of~$E$ over~$K$, then
%$\frac{L'(E/K,1)}{\Omega(E/K)} = \frac{L'(E,1)}{\Omega_E} \cdot
%\frac{L(E_D,1)}{\Omega_{E_D}}$. 
%Also, note that by~\cite[Lem.~2.1]{prasanna:padic}, then considering that $q$ does not divide $c_p(E)$ for any prime~$p$
%that divides~$N$, the same holds for~$c_p(E_D)$. 
%Hence we see that under our hypothesis that
%${rm ord}_q L(E_D,s)/\Omega_{E_D} = 0$, $q$ divides 
%the Birch and Swinnerton-Dyer conjectural
%order of~$\Sha(E/\Q)$.
\end{proof}

In view of Proposition~\ref{prop:agst},
Theorem~\ref{thm:main} and Corollary~\ref{cor:main}
provide theoretical evidence towards 
the Birch and Swinnerton-Dyer conjectural formula~(\ref{gzformula}).
We remark that 
Theorem~\ref{thm:main} is to be compared to
part~(i) of Proposition~\ref{prop:agst}
and Corollary~\ref{cor:main} to part~(ii) of Proposition~\ref{prop:agst}.
%Notice the similarity between the hypotheses 
%Note that the hypotheses of Corollary~\ref{cor:main},
%which shows that certain primes divide the 
%left hand side of
%the Birch and Swinnerton-Dyer conjectural formula~(\ref{gzformula}),
%are somewhat similar (though not identical)
%to the hypotheses of Proposition~\ref{prop:agst},
%which shows that certain primes divide the 
%right hand side of
%conjectural formula~(\ref{gzformula}). 
Regarding the hypothesis in Corollary~\ref{cor:main}
that $q$ does not divide $|J_0(N)(K)_{\rm tor}|$,
we do not know
of any results that would give some criteria on~$q$ which would
imply that this hypothesis holds (unlike the similar situation over~$\Q$,
where at least for prime~$N$, we know by~\cite[Thm~(1)]{mazur:eisenstein}
that $|J_0(N)(\Q)_{\rm tor}|$ equals the numerator of~$\frac{N-1}{12}$).
As in Theorem~\ref{thm:main}, we could have replaced this hypothesis
by the requirement that 
$p$ is coprime to the order of the torsion subgroup of
the projection of~$\T P$ in~$J_0(N)/(I_f \cap I_g) J_0(N)$.
Note that there is some similarity between these hypotheses
and the hypothesis in 
Proposition~\ref{prop:agst} that 
$q$ does not divide 
$|(J_0(N)/\Fdual)(K)_{\rm tor}|$.
In any case, our discussion just above emphasizes
the need to study the torsion in~$J_0(N)$ 
and its quotients over
number fields other than~$\Q$.

The first step in the proof of Theorem~\ref{thm:main} is the following
proposition (also proved in Section~\ref{sec:proofs}), 
which may be of independent interest:
\begin{prop}\label{prop:fact}
Let $f$ be an eigenform of weight~$2$ on~$\Gamma_0(N)$ of analytic rank one
(i.e., whose $L$-function vanishes to order one at~$s=1$).
Let $J'$ be a quotient of~$J=J_0(N)$ through which the map~$J \ra A_f$ factors.
Let $\pi'$ denote the map $J' \ra A_f$ and $\pi''$ the map $J \ra J'$ in
this factorization. Let $F'$ denote the kernel of~$\pi'$. 
Thus we have the following diagram:
$$\xymatrix{
 &  & J \ar[d]_{\pi''} \ar[rd]^{\pi} & & \\
0 \ar[r] & F' \ar[r] & J' \ar[r]^{\pi'}  & A_f \ar[r] & 0\\
}$$
Then
\begin{eqnarray} \label{eqn:fact}
|A_f(K)/ \pi(\T P)| 
 =  \bigg|\frac{J'(K)}{F'(K) + \pi''(\T P)}\bigg| \cdot 
|{\rm ker}\big(H^1(K,F') \ra H^1(K,J')\big)|.
\end{eqnarray}
\end{prop}

\section{Proofs} \label{sec:proofs}
\comment{
Suppose $f$ is congruent modulo a prime ideal~$\pp$ lying over~$p$ 
to another eigenform~$g\in S_2(\Gamma_0(N),\C)$ of
analytic rank~$\geq 3$.
%If $g$ is not new, then the component groups
%interfere (see examples above); 
%so for simplicity, suppose $g$ is new.
Then by ``visibility arguments'', assuming the first part 
of the BSD conjecture, $p$ divides the order of~$\Sha(E/K)$
or some component group. This 
should be reflected on the left hand side of the BSD
conjectural formula~(\ref{gzformula}), and I believe can be proved.
}

In this section, we prove Proposition~\ref{prop:fact}
and Theorem~\ref{thm:main}.
We start by proving Proposition~\ref{prop:fact}.

\begin{proof}[Proof of Proposition~\ref{prop:fact}]
Consider the exact sequence
$0 \ra F' \ra J' \ra A_f \ra 0$. Part of the associated long exact sequence 
of Galois cohomology is
\begin{eqnarray} \label{eqn:les12}
0 \ra F'(K) \ra J'(K) \stackrel{\pi'}{\ra} A_f(K) 
\stackrel{\delta}{\ra} H^1(K, F') \ra H^1(K,J') \ra \cdots,
\end{eqnarray}
where $\delta$ denotes the boundary map.
Note that in this proof, the letters~$\pi'$ and~$\pi''$ denote
$\pi'$ and~$\pi''$ restricted to the $K$-valued points in their domain.
Since $\pi''(\T P) \subseteq J'(K)$, by the exactness of~(\ref{eqn:les12})
we see that 
 $\delta (\pi'(\pi''(\T P)))=0$. Using the exactness of~(\ref{eqn:les12})
again, we see that
$\delta$ thus induces a surjection
$$\phi: A_f(K)/ \pi'(\pi''(\T P)) \ra {\rm ker}\big(H^1(K,F') \ra H^1(K,J')\big).$$
%, and we claim that its
%kernel is isomorphic to $\frac{J(K)}{B(K) + \T P}$. 
%By the exactness of~(\ref{eqn:les1}), 
%and the kernel of~$\phi$ is the image of~$\pi'$ in~$A_f(K)/ \pi'(\pi''(\T P))$.

Since $\pi'(J'(K)) \subseteq {\rm ker}(\delta)$, we see that
$\pi'$ induces a natural map~$\psi: J'(K) \ra \ker(\phi)$.  \\

\noi{\em Claim:} $\psi$ is surjective and its
kernel is $F'(K) + \pi''(\T P)$. 
\begin{proof}
%Note that $\Z \pi(P)
% = (\T/I_f) \pi(P) = \pi(\T P)$. 
Let $x \in J'(K)$. Then 
\begin{tabbing}
$x \in \ker(\psi)$ \= 
$\iff \pi'(x) = 0 \in \ker(\phi) \hookrightarrow A_f(K)/ \pi'(\pi''(\T P))$\\
\> $\iff \pi'(x) \in \pi' (\pi''(\T P))$\\
\> $\iff \exists t \in \T : x - \pi''(t P) \in \ker(\pi') = F'(K)$\\
\> $\iff x \in F'(K) + \pi''(\T P)$. 
\end{tabbing}
Thus $\ker(\psi) = F'(K) + \pi''(\T P)$.
To prove surjectivity of~$\psi$, note that given an
element of~$\ker(\phi)$, we can write the element
 as $y + \pi'(\pi''(\T P))$
for some $y \in A_f(K)$ such that $\delta(y) = 0$.
Then by the exactness of~(\ref{eqn:les12}), 
$y \in {\rm Im}(\pi')$, hence 
$y + \pi'(\pi''(\T P)) \in {\rm Im}(\psi)$.
\end{proof}

By the discussion above,  we get an exact sequence:
\begin{eqnarray} \label{eqn:finite}
0 \ra \frac{J'(K)}{F'(K) + \pi''(\T P)} \stackrel{\psi'}{\ra}
\frac{A_f(K)}{\pi'(\pi''(\T P))} \stackrel{\phi}{\ra} 
{\rm ker}\big(H^1(K,F') \ra H^1(K,J')\big) \ra 0,
\end{eqnarray}
where $\psi'$ is the natural map induced by~$\psi$.
Now
$$|A_f(K)/ \pi'(\pi''(\T P))| = |A_f(K)/ \pi(\T P)|\ ,$$
and the latter group is finite in our situation.
Hence all groups in~(\ref{eqn:finite}) are finite, and 
Proposition~\ref{prop:fact} now follows from the exactness
of~(\ref{eqn:finite}).
\end{proof}

In the rest of this section, we prove Theorem~\ref{thm:main}.
We work in slightly more generality in the beginning
and assume that $f$ and~$g$ are any newforms 
(whose Fourier coefficients need not be integers). 
Thus the associated newform quotients~$A_f$ and~$A_g$ need not be elliptic
curves, but we will still denote them by~$E$ and~$F$ (respectively)
for simplicity of notation.
%Recall that $E$ and~$F$ are optimal elliptic curves with associated
%newforms~$f$ and~$g$ such that

Recall that $I_g = {\rm Ann}_{\scriptscriptstyle \T} g$, and 
let $J' = J/(I_f \cap I_g)J$. The quotient map $\pi: J \ra J/I_f J$
factors through~$J'= J/(I_f \cap I_g)J$. 
By proposition~(\ref{prop:fact}), it suffices to show 
that $r$ divides 
%one of the two factors on 
the right side of~(\ref{eqn:fact}) with our choice of~$J'$, 
which is what we will do. 

Recall that $\pi'$ denotes 
the projection $J' \ra A_f = E$,
$F' = \ker \pi'$, and $\pi''$ denotes the projection $J \ra J'$. 
Let $B$ denote the kernel of~$\pi: J \ra E$, which is the
abelian subvariety~$I_f J$ of~$J$.
We have the following diagram, in which the two sequences of four arrows are 
exact (one horizontal and one upwards diagonal):

$$\xymatrix{
 & \Fdual \ar@{^(->}[rd] & \quad \Edual \ar@{^(->}[d] \ar[rd]^{\sim} & & 0 \\
0 \ar[r] & B \ar[dd] \ar[r] & J \ar[r]^{\pi} \ar[d]_{\pi''} & E \ar[r] 
\ar[ur] & 0\\
& & J'\ar[rd] \ar[ur]_{\pi'} &  & \\
& F' \ar[ur] & & F  & \\
0 \ar[ur] & & & & \\
}$$

Now $F'$ is connected, since it is
a quotient of~$B$ (as a simple diagram chase above shows) and $B$ is connected. 
Thus, by looking at dimensions, one sees that $F'$ is the image of~$\Fdual$
under~$\pi''$. Since the composite $\Fdual \hookrightarrow J \ra J' \ra F$
is an isogeny, the 
the quotient map $J' \ra F$ induces an isogeny
$\pi''(\Fdual) \sim F$, and hence 
an isogeny $F' \sim F$. Thus $F'$ and~$F$ have the same rank 
(over~$\Q$ or over~$K$).
Let $E'$ denote $\pi''(\Edual)$. Since $\pi$ induces an isogeny
from~$\Edual$ to~$E$, we see that $\pi'$ also induces an isogeny
from~$E'$ to~$E$. Thus $E'$ and~$E$ have the same rank 
(over~$\Q$ or over~$K$).

Now we impose the assumption that $f$ and~$g$ have integer
Fourier coefficients, so that $E$ and~$F$ are elliptic curves.
Recall that $r$ is an integer such that
(a) $f$ and~$g$  are congruent modulo~$r$ and 
for every prime~$p$ that divides~$r$, $f$ and~$g$ are not
congruent modulo a power of~$p$ greater than~$p^{{\rm ord}_p r}$, 
(b) all maximal ideals~$\m$ of~$\T$ whose residue characteristic
divides~$r$ and that are in the support of~$J[I_f]$
or of~$J[I_g]$ satisfy multiplicity one, and
(c) $r$ is coprime to the order of the torsion subgroup of
the projection of~$\T P$ in~$J/(I_f \cap I_g) J$.
%Then $r$ divides $|E(K)/ \Z \pi(P)|$.
Our goal is to show that 
$r$ divides $|E(K)/ \Z \pi(P)|$.

From Lemma~\ref{lem:containment} and on applying~$\pi''$, we see that 
$E'[r] = F'[r]$ and both are
direct summands 
of~$E' \cap F'$ as~${\rm Gal}(\Qbar/\Q)$-modules.
In particular, the 
natural maps $H^1(K, E'[r]) \rightarrow H^1(K, E' \cap F')$
and $H^1(K, F'[r]) \rightarrow H^1(K, E' \cap F')$ are
injections. 
Recall that $E$ has analytic rank one and $F$ has Mordell-Weil rank more
than one. Then the abelian group~$F(K)$ has rank more than one, and 
as remarked just before Conjecture~\ref{gzformula}, the abelian group~$E(K)$ 
has rank one. Also, note that the newform~$g$ has analytic rank greater 
than one, since otherwise the Mordell-Weil
rank of~$F$ would be at most one. 
With an eye towards potential generalizations,
we remark that after this paragraph, we will not explicitly
use the fact that $E$ and $F$ have dimension one (i.e., are elliptic
curves). Thus if the conclusions of this paragraph are satisfied,
then the rest of the argument would go through even if $E$ and~$F$ have
dimension greater than one.

Consider the following commutative diagram, where
the top and bottom rows are the Kummer
exact sequences of~$E'$ and~$F'$ respectively, 
and the other maps are the obvious natural maps:

$$\xymatrix{
0  \ar[r] & E'(K)/ r E'(K) \ar[r] &  H^1(K, E'[r]) 
 \ar[r] \ar@{^(->}[rdd] &  H^1(K, E')[r] \ar@{^(->}[rd] \ar[r] &  0\\
& & & &  H^1(K, E') & \\
 & & & H^1(K, E' \cap F') \ar[ur] \ar[dr]&   & \\
& & & &  H^1(K, F') & \\
0  \ar[r] & F'(K)/ r F'(K) \ar[r]^{\delta'} &  H^1(K, F'[r]) \ar@{=}[uuuu]
 \ar[r] \ar@{^(->}[ruu] &  H^1(K, F')[r] \ar@{^(->}[ru] \ar[r] &  0\\
}$$

Let $Q$ be a generator for the free part of~$E'(K)$
(which is isomorphic to the free part of~$E(K)$). Then
from the top exact sequence in the diagram above,
we see that $Q$ gives rise
to a non-trivial element~$\sigma$ in~$H^1(K, E'[r])$. 
%and hence a non-trivial element~$\sigma$ in~$H^1(K, E'[p])$.
\comment{
Similarly a point of infinite order in~$F'(K)$ gives rise
to a non-trivial element in~$H^1(K, F'[r])$; 
let $S$ be the set of the images of these points in
$H^1(K, F'[r]) = H^1(K, E'[r]) = H^1(K, F'[r]) $.

We consider two cases corresponding to whether
$\sigma$ is in the subset $\delta'(F'(K)/r F'(K))$
of~$H^1(K, F'[r]) = H^1(K, E'[r])$  or not
(where $\delta'$ is the boundary map in the Kummer exact sequence
associated to~$F'$, as indicated in the diagram above).\\

\noindent Case 1: $\sigma \in \delta'(F'(K)/r F'(K))$ \\
}

Let $r'$ be the smallest positive integer such that 
$r' \sigma \in \delta'(F'(K)/r F'(K))$ 
(where $\delta'$ is the boundary map in the Kummer exact sequence
associated to~$F'$, as indicated in the diagram above).
Thus $r'$ divides~$r$ (since $r \sigma = 0 \in \delta'(F'(K)/r F'(K))$).
Then, by the top and bottom exact sequences in the diagram above, 
$r' \sigma$ 
maps to the trivial element in both
$H^1(K,E')[r]$ and $H^1(K,F')[r]$, and hence in
$H^1(K,E')$ and $H^1(K,F')$. The image of $r' \sigma$ in~$H^1(K,E' \cap F')$
is then a non-trivial element of order~$r/r'$
that dies in~$H^1(K, F')$ and in~$H^1(K, E')$.
Thus we see that $r/r'$ divides the order of
$\ker(H^1(K,E' \cap F') \ra H^1(K,E' \oplus F'))$.

\begin{lem} We have
$$
\frac{J'(K)}{F'(K) + E'(K)} \isom 
\ker\big(H^1(K,E' \cap F') \ra H^1(K,E' \oplus F')\big).
$$
\end{lem}

\begin{proof}
Following a similar situation in~\cite{cremona-mazur}, 
consider the short exact sequence
\begin{eqnarray} \label{eqn:ses}
0 \ra E' \cap F' \ra E' \oplus F' \ra J' \ra 0, 
\end{eqnarray}
where the map $E' \cap F' \ra E' \oplus F'$
is the anti-diagonal embedding $x \mapsto (-x,x)$
and the map $E' \oplus F' \ra J'$ is given by $(x,y) \mapsto x+y$.
Part of the associated long exact sequence is
$$\cdots \ra E'(K) \oplus F'(K) \ra J'(K) \ra 
H^1(K,E' \cap F') \ra H^1(K,E' \oplus F') \ra \cdots,$$
from which we get the lemma.
\end{proof}

By the lemma and the discussion preceding it, 
we see that
%From the previous paragraph,
%$\sigma$ is a non-trivial element of order~$r$ of
%the right side of~(\ref{galcoh2}). Hence we see
%that 
$r/r'$ divides~$|\frac{J'(K)}{F'(K) + E'(K)}|$. 
%In order to prove Theorem~\ref{thm:main}, we may assume that $r$ does
%not divide~$E(K)_{\rm tor}$, since if it does, then considering
%that $\pi(P)$ has infinite order, $r$ divides
%$|E(K)/ \Z \pi(P)|$.
If $h$ is an eigenform of weight~$2$ on~$\Gamma_0(N)$,
then $\T P \cap \Ahdual(K)$ is infinite if and only if $h$ has
analytic rank one (this follows by~\cite[Thm~6.3]{gross-zagier}
if $h$ has analytic rank bigger than one,
and the fact that $\Ahdual(K)$
is finite if~$h$ has analytic rank one, by~\cite{kollog:finiteness}).
Thus, considering that $J'$ is isogenous to $E' \oplus F'$ and
$g$ has analytic rank greater than one, we
see that the free parts of $E'(K)$ and of~$\pi''(\T P)$ agree. 
Since we are assuming that 
$r$ is coprime to the order of the torsion part of~$\pi''(\T P)$
(which is condition~(c) on~$r$),
considering that $r/r'$ divides
$|\frac{J'(K)}{F'(K) + E'(K)}|$,
we see that $r/r'$ divides
$|\frac{J'(K)}{F'(K) + \pi''(\T P)}|$ as well.
Thus $r/r'$ divides the
the first factor on the right side of~(\ref{eqn:fact}).

%\noindent Case 2: $\sigma \not\in \delta'(F'(K)/r F'(K))$ \\
If $r'=1$, then we are done, so let us assume that $r' > 1$.
Then $\sigma \not\in \delta'(F'(K)/r F'(K))$. 
So while the image of~$\sigma$ in~$H^1(K,E')[r]$ is trivial by 
the top exact sequence in the diagram above, 
the image of~$\sigma$ in~$H^1(K,F')[r]$ generates a subgroup of order~$r'$,
%in~$H^1(K,F')[r]$%, hence in~$H^1(K,F')$, 
by the lower exact sequence of the diagram above
(recall that $r'$ is the {\em smallest} positive integer such that 
$r' \sigma \in \delta'(F'(K)/r F'(K))$).
%by the Kummer exact sequence
%of~$F'$, considering that $\sigma$ does not arise from~$F'(K)$,
%keeping in mind that $r$ is coprime to the order of~$F'(K)_{\rm tor}$. 
Thus, from the same diagram, we see that 
$\sigma$, viewed as an element of~$H^1(K,E' \cap F')$,
maps to the trivial element of~$H^1(K,E')$ 
but a nontrivial element~$\sigma'$ of~$H^1(K,F')$ 
of order~$r'$.
Following a similar situation in~\cite{mazur:visthree},
considering the exactness of
$$H^1(K, E' \cap F') \ra H^1(K,E') \oplus H^1(K,F') \ra 
H^1(K,J'),$$
which is part of the long exact sequence associated to~(\ref{eqn:ses}),
we see that the element $(0,\sigma')$ in the middle
group dies in the rightmost group, since it arises from the
element~$\sigma$ in the leftmost group.
Thus $\sigma' \in H^1(K,F')$ dies in~$H^1(K,J')$, and thus is a nontrivial 
element of order~$r'$ of~${\rm ker}\big(H^1(K,F') \ra H^1(K,J')\big)$.
Hence $r'$ divides the second factor on the right side of~(\ref{eqn:fact}). 

%In either of the two cases above, we see that $r$ divides
%one of the two factors on
Thus $r/r'$ divides the first factor on the right side of~(\ref{eqn:fact}), 
and $r'$ divides the second factor on the right side of~(\ref{eqn:fact}),
and so $r$ divides the right side of~(\ref{eqn:fact}), as was to be shown.

\bibliographystyle{amsalpha}         

\providecommand{\bysame}{\leavevmode\hbox to3em{\hrulefill}\thinspace}
\providecommand{\MR}{\relax\ifhmode\unskip\space\fi MR }
\providecommand{\MRhref}[2]{%
  \href{http://www.ams.org/mathscinet-getitem?mr=#1}{#2}
}
\providecommand{\href}[2]{#2}

\end{document}